\documentclass[10pt]{amsart}
\usepackage{amsthm,amsmath,amssymb}

\global\overfullrule10pt

\theoremstyle{plain}
\newtheorem{theorem}{Theorem}
\newtheorem{corollary}[theorem]{Corollary}
\newtheorem{lemma}[theorem]{Lemma}
\newtheorem{proposition}[theorem]{Proposition}
\newtheorem*{theorem*}{Theorem}
\newtheorem*{corollary*}{Corollary}
\theoremstyle{definition}

\theoremstyle{remark}
\newtheorem{remark}[theorem]{Remark}

\newcommand\CC{{\mathbf C}}
\newcommand\RR{{\mathbf R}}
\newcommand\ZZ{{\mathbf Z}}
\newcommand\NN{{\mathbf N}}
\newcommand\DD{{\mathbf D}}
\newcommand\BB{{\mathbf B}}
\newcommand\spr[1]{\langle#1\rangle}
\newcommand\hol{{\text{hol}}}

\newcommand\mh{{\text{Mh}}}
\newcommand\ph{{\text{ph}}}
\renewcommand\Re{\operatorname{Re}}

\newcommand\oz{\overline z}
\newcommand\Bn{{\BB^n}}
\newcommand\pBn{{\partial\Bn}}
\newcommand\Hardy{{\text{\rm Hardy}}}

\newcommand\wt{\widetilde}
\newcommand\td{\wt\Delta}

\newcommand\FF[5]{{}_#1\!F_#2\Big(\begin{matrix}#3\\#4\end{matrix}\Big|#5\Big)}
\newcommand\FE[4]{F_#1\Big(\begin{matrix}#2\\#3\end{matrix}\Big|#4\Big)}
\newcommand\Mh{$M$-harmonic }

\newcommand\GG[2]{\Gamma\Big(\begin{matrix}#1\\#2\end{matrix}\Big)}
\newcommand\chpq{\mathcal H^{pq}}
\newcommand\bhpq{\mathbf H^{pq}}
\newcommand\Un{{U(n)}}

\newcommand\dsph{\Delta_{\text{sph}}}
\newcommand\cR{\mathcal R}
\newcommand\cN{\mathcal N}

\makeatletter
\newcommand{\vast}{\bBigg@{4}}
\newcommand{\Vast}{\bBigg@{5}}
\makeatother

\renewcommand\[{\begin{equation}}
\renewcommand\]{\end{equation}}
\allowdisplaybreaks[4]

\renewcommand\AA{\mathcal A}

\newcommand\autBn{\operatorname{Aut}(\Bn)}
\newcommand\MM{\mathcal M}
\newcommand\dbar{\overline\partial}

\newcommand\oL{\overline L}
\newcommand\nadseb[3]{\Big(\begin{matrix}#1\\#2\end{matrix}\Big|#3\Big)}

\begin{document}

\title[$M$-harmonic Besov spaces]{Besov-Bergman spaces of $M$-harmonic functions}
\author[P.~Blaschke]{Petr Blaschke}
\address{Mathematics Institute, Silesian University in Opava,
 Na~Rybn\'\i\v cku~1, 74601~Opava, Czech Republic}
 \email{Petr.Blaschke{@}math.slu.cz}
\author[M.~Engli\v s]{Miroslav Engli\v s}
\address{Mathematics Institute, Silesian University in Opava,
 Na~Rybn\'\i\v cku~1, 74601~Opava, Czech Republic {\rm and }
 Mathematics Institute, \v Zitn\' a 25, 11567~Prague~1,
 Czech Republic}
\email{englis{@}math.cas.cz}
\thanks{Research supported by GA\v CR grant no.~25-18042S
 and RVO funding for I\v CO~67985840.}
\begin{abstract} We~show that the weighted Bergman spaces of M-harmonic functions 
(functions annihilated by the invariant Laplacian on the unit ball of the complex n-space),
as~well as their analytic continuation (in~the spirit of Rossi and Vergne), coincide with
the certain Besov-type spaces, which were studied by Folland. Characterizations in terms
of tangential derivatives are given, and for appropriate values of the weight parameter,
these spaces are also shown to coincide with the subspaces of all M-harmonic fucntions
in the Sobolev space of order~$t$ on the~ball, $0\le t\le n$. Unlike the holomorphic case,
the~last result is shown to fail in general for other values of~$t$. The~main tool in the 
proofs are asymptotic estimates for certain integrals of squared hypergeometric functions, 
which seem to be of interest in their own right and may find other applications.
\end{abstract}

\maketitle

\section{Introduction}
Let $\Bn$ denote the unit ball in the complex $n$-space~$\CC^n$, $n\ge1$,
with its hyperbolic Laplacian
$$ \td = 4(1-|z|^2) \sum_{j,k=1}^n (\delta_{jk}-z_j\oz_k) \frac{\partial^2}{\partial z_j\partial\oz_k}  $$
invariant under the group $\autBn$ of all biholomorphic self-maps (Moebius transformations) of~$\Bn$.
Functions annihilated by $\td$ are called \emph{Moebius harmonic} (or~\emph{invariantly harmonic}),
or~\emph{\Mh}for short; see Chapter~4.3 in Rudin~\cite{Ru}. For any $s>-1$, one has the standard
rotationally symmetric probability measure on~$\Bn$
$$ d\mu_s(z) = \frac{(s+1)_n}{\pi^n} (1-|z|^2)^n \, dz    $$
(where $dz$ denotes the Lebesgue measure on~$\CC^n$), and the associated subspaces of all \Mh functions
in~$L^2(\Bn,d\mu_s)$, or~\emph{weighted \Mh Bergman spaces}:
\[ \MM_s := \{ f\in L^2(\Bn,d\mu_s): f\text{ is \Mh on }\Bn \} .  \label{TA} \]
Here $(\nu)_k:=\nu(\nu+1)\dots(\nu+k-1)$ stands for the usual Pochhammer symbol (rising factorial).

For the analogous spaces of holomorphic, rather than $M$-harmonic, functions,
i.e.~the ordinary weighted Bergman spaces
$$ \AA_s := \{ f\in L^2(\Bn,d\mu_s): f\text{ is holomorphic on }\Bn \}, \qquad s>-1,  $$
the corresponding norm can easily be expressed in terms of Taylor coefficients:
\[ \|f\|_s^2 := \int_\Bn |f(z)|^2 d\mu_s(z) = \sum_{\alpha\text{ multi-index}}
 |f_\alpha|^2 \frac{\alpha!}{(s+n+1)_{|\alpha|}} \label{TB}  \]
for $f=\sum_\alpha f_\alpha z^\alpha$
(with the usual multi-index notation). Remarkably, the right-hand side continues to be
(well-defined and) positive definite not only for $s>-1$, but for all $s>-n-1$; in~this way,
one~obtains an ``analytic continuation'' of the weighted Bergman spaces $\AA_s$, $s>-1$,
to all $s>-n-1$, with the associated squared norms being still given by the sum in~\eqref{TB}.
This has been done even in much greater generality, for $\Bn$ replaced by an arbitrary
bounded symmetric domain in~$\CC^n$, by~Rossi and Vergne~\cite{RV}.
Furthermore, the~corresponding spaces $\AA_s$, $s>-n-1$, turn out to be Besov-type spaces:
for~our case of~$\Bn$, $\AA_s$~coincides with the subspace $W^{-s/2}_\hol(\Bn)$ of all
holomorphic functions in the Sobolev space $W^{-s/2}(\Bn)$ on $\Bn$ of order~$-\frac s2$,
with equivalent norms. This makes it possible to extend the definition of $\AA_s$ even further,
from $s>-n-1$ to all real~$s$.
All~this can actually be shown to remain in force also for $\Bn$ replaced by an arbitrary
bounded strictly-pseudoconvex domain $\Omega$ in~$\CC^n$ with smooth boundary.
The~situation also turns out to be completely analogous for the weighted Bergman spaces of
\emph{harmonic} functions on a bounded domain in the real $n$-space $\RR^n$ with smooth boundary
(we~omit the details).

For~our \Mh weighted Bergman spaces $\MM_s$ above, it~has recently been shown by E.-H.~Youssfi
and one of the authors \cite{EY2} that there~is, again, an~analogous ``analytic continuation''
of the spaces $\MM_s$ from $s>-1$ to $s>-n-1$; furthermore, as a ``residue'' at $s=-n-1$,
one~obtains a certain Dirichlet space of \Mh functions. Our~goal in this paper is to show
that these ``analytically continued'' spaces $\MM_s$, $s>-n-1$, are again Besov (or~Sobolev)
type spaces, however, with derivatives taken only in complex tangential directions.
Furthermore, their definition can again be extended even to all real~$s$, and then they
coincide also with the subspaces $W^{-s/2}_\mh(\Bn)$ of all \Mh functions in the ordinary
Sobolev space $W^{-s/2}(\Bn)$ on $\Bn$ again of order $-\frac s2$, for $-2n\le s\le0$,
but in general this is no longer true for other values of~$s$.

To~describe our results, consider the vector fields
\[ L_{jk} := \oz_j\frac\partial{\partial z_k} - \oz_k\frac\partial{\partial z_j}, \quad
 \oL_{jk} := z_j\frac\partial{\partial\oz_k} - z_k\frac\partial{\partial\oz_j},
 \qquad j,k=1,\dots,n,  \label{TC} \]
on~$\CC^n$; clearly, these are tangential to all spheres $|z|\equiv$const., and $L_{jk}$
and $\oL_{jk}$ generate (very redundantly) the holomorphic and the anti-holomorphic
complex tangent space, respectively, to~these spheres. For~notational convenience,~set
\[ L_{jk}=:X_{(n-1)j+k}, \quad \oL_{jk}=:X_{n^2+(n-1)j+k}, \qquad j,k=1,\dots,n,  \label{TD} \]
so~that $\{L_{jk}\}_{j,k=1}^n\cup\{\oL_{jk}\}_{j,k=1}^n=\{X_j\}_{j=1}^{2n^2}$. The~operator
$$ \square := -\sum_{j,k=1}^n (L_{jk}\oL_{jk}+\oL_{jk}L_{jk})  $$
can~be viewed as the counterpart on $\pBn$ of the Folland-Stein sublaplacian on
the Heisenberg group~\cite{FS}. Although the individual operators $L_{jk}$ and $\oL_{jk}$
do~not necessarily map an \Mh function again into an \Mh function, the operator~$\square$ does;
in~fact, one~can exhibit an explicit orthonormal basis of eigenfunctions of $\square$
with nonnegative eigenvalues, thanks to which it is possible to define $(I+\square)^t f$
for any $t\in\CC$ and any \Mh function $f$ on~$\Bn$. (See Section~\ref{Sec2} below for all details.)
Finally, in~addition to the norms $\|\cdot\|_s$, $s>-1$, in~$L^2(\Bn,d\mu_s)$, let~us also denote
by $\|\cdot\|_\pBn$ the norm in $L^2(\pBn,d\sigma)$ with respect to the normalized surface
measure $d\sigma$ on the unit sphere~$\pBn$, and for an \Mh function~$f$, denote
\[ \|f\|^2_\Hardy := \sup_{0<r<1} \|f_r\|^2_\pBn, \quad\text{where } f_r(\zeta):=f(r\zeta),
 \quad 0<r<1, \;\zeta\in\pBn.  \label{TE}  \]
Our~main result is the following.

\begin{theorem*}[Theorem~\ref{PA}]
For $f$ \Mh on $\Bn$, $n\ge2$, and $s>-n-1$, the following assertions are equivalent.
\begin{itemize}
\item[(a)] $f\in\MM_s$;
\item[(b)] $\sum_{j_1,\dots,j_m=1}^{2n^2} \|X_{j_1}\dots X_{j_m}f\|^2_{s+m} + |f(0)|^2 <+\infty$
for some (equivalently, any) $m=0,1,2,\dots$, $m>-s-1$;
\item[(c)] $\|(I+\square)^{t/2}f\|^2_{s+t}<+\infty$ for some (equivalently, any) $t>-s-1$;
\item[(d)] $\|(I+\square)^{-(s+1)/2}f\|^2_\Hardy<+\infty$.
\end{itemize}
Furthermore, the quantities above are all equivalent to $\|f\|^2_s$.
\end{theorem*}

In~particular, defining, for any $s\in\RR$,
$$ \wt\MM_s := \{ f\text{ \Mh on }\Bn: \;\|f\|_{\tilde s}:=\|(I+\square)^{-(s+1)/2}f\|_\Hardy<+\infty \}, $$
we~will have $\wt\MM_s=\MM_s$ for all $s>-n-1$, with equivalent norms, yielding again an extension
of the definition of $\MM_s$ from $s>-n-1$ to all real~$s$.
Furthermore, the~quantities in parts (b), (c) of Theorem~\ref{PA} will still be equivalent norms on~$\wt\MM_s$.

We~pause to note that the spaces $\wt\MM_s$ were studied in the paper by Folland~\cite{Fo1},
where they were denoted by~$B_{-s-1,0}$ (cf.~page~111 there); thus our result implies that
for $s>-1$ these Besov spaces of Folland actually coincide with the weighted \Mh Bergman spaces~$\MM_s$.

Our~second result is a characterization of the ``full'' Sobolev spaces of \Mh functions
$$ W^t_\mh(\Bn) := \{f\in W^t(\Bn): \;f\text{ is \Mh on }\Bn \},  $$
albeit only for a certain range of the real order~$t$.

\begin{theorem*}[Theorem~\ref{PB}]
For $0\le t\le n$, $W^t_\mh(\Bn) = \wt\MM_{-2t}$.
\end{theorem*}

The~main point of the last theorem~is, of~course, that --- perhaps surprisingly ---
the~complex tangential derivatives $L_{jk},\oL_{jk}$ already control (on~\Mh functions)
also the complex normal derivative and the real normal derivative.
Naturally, this fails for general (i.e.~not $M$-harmonic) functions.

The~assertion of Theorem~\ref{PB} is also shown to fail in general for $t>n$.

Our~proofs use the decomposition of \Mh functions into ``bigraded spherical harmonics''
due to Folland~\cite{Fo1},~\cite{Fo2}. Namely, any \Mh function on $\Bn$ can be uniquely
written in the form
\[ f = \sum_{p,q=0}^\infty f_{pq} ,   \label{TF}  \]
where the series converges uniformly on compact subsets, and the ``pieces'' $f_{pq}$
have the special form
$$ f_{pq}(z) = S_{pq}(|z|^2) \tilde f_{pq}(z,\oz)  $$
where $S_{pq}$ is a certain hypergeometric function (independent of~$f$) while $\tilde f_{pq}$
is~a (Euclidean) harmonic polynomial on $\CC^n$ homogeneous of bidegree $(p,q)$ in~$(z,\oz)$.
Furthermore, for~$f$ as in~\eqref{TF}, the~norm in~$\MM_s$, $s>-1$, is~given~by
\[ \|f\|^2_s = \sum_{p,q=0}^\infty c_{pq}(s) \|\tilde f_{pq}\|^2_\pBn  \label{TH}  \]
where
\[ c_{pq}(s) = \frac{(s+1)_n}{\Gamma(n)} \int_0^1 t^{p+q+n-1} (1-t)^s S_{pq}(t)^2 \,dt. \label{TI} \]
It~was shown in \cite{EY2} that $c_{pq}(s)$ actually extend, for each $p,q\ge0$, to~a~holomorphic
function of $s$ on the entire~$\CC$, except for (possible) poles at $s=-n-1-j$, $j=0,1,2,\dots$.
It~turns out to be surprisingly difficult to obtain bounds for the last integral with $s>-1$,
valid uniformly for all $p,q\ge0$, and even more so for its analytic continuation to $s>-n-1$
just mentioned. The~following result is the crucial ingredient for the proof of Theorem~\ref{PA}.

\begin{corollary*}[Corollary~\ref{PC}]
For each $s>-n-1$, there exist constants $0<c_s<C_s<+\infty$ such that
\[ c_s \le (p+1)^{s+1}(q+1)^{s+1} c_{pq}(s) \le C_s \qquad\forall p,q\ge0.  \label{TG} \]
\end{corollary*}

This estimate also resolves (affirmatively) a conjecture from the earlier paper~\cite{EY2}.
The~corollary is an easy consequence of an alternative formula for~$c_{pq}(s)$ (Theorem~\ref{PE}),
which is of interest in its own right and may find other applications.

The~proof of Corollary~\ref{PC} is the content of Section~\ref{Sec3}, after recalling the needed
preliminaries on \Mh functions in Section~\ref{Sec2}. The~proofs of Theorems~\ref{PA} and~\ref{PB}
appear in Sections~\ref{Sec4} and~\ref{Sec5}, respectively.

Throughout the paper, we~use the notation $\spr{f,g}$ for the inner product in various function spaces.
The notation $A\lesssim B$
means that there exists a finite constant~$c$, independent of the variables in question, such that $A\le cB$;
and $A\asymp B$ means that $A\lesssim B$ and $B\lesssim A$. The~symbols $\frac\partial{\partial z_j}$
and $\frac\partial{\partial\oz_j}$, commonly abbreviated just to $\partial_j$ and $\dbar_j$, respectively,
stand for the usual Wirtinger operators on~$\CC^n$. For~notational convenience, we~sometimes use the shorthand
$$ \GG{a_1,a_2,\dots,a_j}{b_1,b_2,\dots,b_k} := \frac{\Gamma(a_1)\Gamma(a_2)\dots\Gamma(a_j)}{\Gamma(b_1)\Gamma(b_2)\dots\Gamma(b_k)}. $$
Finally, $\ZZ,\NN,\RR$ and $\CC$ denote the sets of all integers, all nonnegative integers, all real
and all complex numbers, respectively.

\section{Preliminaries} \label{Sec2}
Composition with elements of the unitary group $\Un$
\[ f\longmapsto f\circ U^{-1}, \qquad U\in\Un,  \label{UE}  \]
gives a unitary representation of $\Un$ on the space $L^2(\pBn,d\sigma)$.
The~Peter-Weyl decomposition of this action into irreducible components is given~by
\[ L^2(\pBn,d\sigma) = \bigoplus_{p,q=0}^\infty \chpq,  \label{UA}  \]
where~$\chpq$, the spaces of ``bigraded spherical harmonics'', are defined~as
$$ \chpq := \{\tilde f|_\pBn: \; \tilde f\in\wt\chpq \},  $$
where $\wt\chpq$ is the vector space of all harmonic polynomials $\tilde f(z,\oz)$
on $\CC^n$ homogeneous of degree $p$ in $z$ and homogeneous of degree $q$ in~$\oz$.
The~restriction map $\tilde f\to f=\tilde f|_\pBn$ is one-to-one from $\wt\chpq$ onto~$\chpq$.

Denote
\[ S_{pq}(t) := \frac{\FF21{p,q}{p+q+n}t}{\FF21{p,q}{p+q+n}1} =
 \frac{(n)_p(n)_q}{(n)_{p+q}} \FF21{p,q}{p+q+n}t , \label{UB} \]
where ${}_2\!F_1$ denotes the Gauss hypergeometric function; and define the space of
``solid harmonics'' of bidegree $(p,q)$~by
\[ \bhpq := \{S_{pq}(|z|^2) \tilde f(z,\oz): \; \tilde f\in\wt\chpq \}.  \label{UC}  \]
Then $S_{pq}(1)=1$ for all $p,q\ge0$, and each $f_{pq}\in\bhpq$ is \Mh on $\Bn$ and
coincides with $\tilde f(z,\oz)$ on~$\pBn$.
Furthermore, any \Mh function $f$ on $\Bn$ can be uniquely decomposed~as
\[ f = \sum_{p,q=0}^\infty f_{pq}, \qquad f_{pq}\in\bhpq,  \label{UD}  \]
with convergence uniform on compact subsets of~$\Bn$.

\begin{remark} Loosely speaking, a~harmonic function $\sum_{p,q}\tilde f_{pq}(z,\oz)$,
$\tilde f_{pq}\in\wt\chpq$, is turned into an \Mh function $\sum_{p,q} S_{pq}(|z|^2) \tilde f_{pq}(z,\oz)$
by multiplying it by $S_{pq}(|z|^2)$ on each~$\wt\chpq$. \qed
\end{remark}

We~refer the reader to Rudin~\cite[Sections~12.1--12.2]{Ru}, Krantz~\cite[Sections~6.6--6.8]{Krn}
and Ahern, Bruna and Cascante \cite{ABC} for the material above; basic ingredients go back to
Folland~\cite{Fo1} \cite{Fo2}. A~slight caveat about notation: for~later convenience, we~are using
the notation $S_{pq}$ (see~\eqref{UB} above) instead of $S^{pq}(r):=r^{p+q}S_{pq}(r^2)$ in the above references.

Recall that in polar coordinates $z=r\zeta$ on~$\CC^n$ ($r>0$, $\zeta\in\pBn$),
the~Euclidean Laplacian $\Delta$ is given~by
\[ \Delta = \frac{\partial^2}{\partial r^2} + \frac{2n-1}r \frac\partial{\partial r} + \frac1{r^2}\dsph, \label{UL} \]
where $\dsph$ is the \emph{spherical Laplacian}, which involves only differentiations
with respect to the $\zeta$ variables. The~operator $\dsph$ can be expressed explicitly~as
$$ \dsph = - \cR^2 + \sum_{j,k=1}^n (L_{jk}\oL_{jk}+\oL_{jk}L_{jk}) = -\cR^2 - \square  $$
where $\cR$ stands for the complex normal derivative
\[ \cR := \sum_{j=1}^n \Big(z_j \frac\partial{\partial z_j} - \oz_j \frac\partial{\partial\oz_j}\Big) . \label{UR}  \]
Both $\dsph$ and $\cR$ --- and, hence, also~$\square$ --- commute with the action \eqref{UE} of~$\Un$.
From the irreducibility of the multiplicity-free decomposition \eqref{UA} it therefore
follows by abstract theory (Schur lemma) that $\dsph$, $\cR$ and $\square$ map each~$\chpq$
(and~$\bhpq$) into itself and actually reduce on it to a multiple of the identity.
Evaluation on e.g.~the element $\zeta_1^p\overline\zeta_2^q\in\chpq$ (for $n\ge2$) shows that, explicitly,
\[ \begin{aligned}
 \dsph | \chpq &= - (p+q)(p+q+2n-2) I | \chpq, \\
 \cR | \chpq &= (p-q)I | \chpq, \\
 \square | \chpq &= (4pq+(2n-2)(p+q)) I | \chpq , \end{aligned}  \label{UF}  \]
and the same holds, thanks to~\eqref{UC}, for $\bhpq$ in the place of~$\chpq$.
(These formulas also prevail for $n=1$; in~that case, $\dsph=-\cR^2$, $\square=0$,
and $\chpq=\{0\}$ unless $pq=0$, while $\mathcal H^{p0}=\CC z^p$ and $\mathcal H^{0q}=\CC\oz^q$.)

One~consequence of \eqref{UF} is~that, in~particular, $\dsph$, $\cR$ and $\square$ commute.
The~analogue of \eqref{UL} for the invariant Laplacian,
$$ \td = (1-r^2)^2 \frac{\partial^2}{\partial r^2} + \Big(\frac{2n-1}r-r\Big) \frac\partial{\partial r} + \Big( \frac{\dsph}{r^2}+\cR^2\Big)I , $$
together with the fact that $\dsph$, $\cR$ and $\square$ involve only tangential differentiations,
thus shows that these three operators commute with~$\td$. In~particular, they preserve the kernel of~$\td$,
i.e.~$\dsph$, $\cR$ and $\square$ map \Mh functions again into \Mh functions.
(Of~course, this also follows directly from \eqref{UD} and~\eqref{UF}.)

Using the formulas~\eqref{UF}, one can define $(I+\square)^t f$ for any $t\in\CC$ and $f$ \Mh on~$\Bn$ by setting
\[ (I+\square)^t \sum_{p,q} f_{pq} := \sum_{p,q} [4pq+(2n-2)(p+q)+1]^t f_{pq} , \label{UG} \]
for $f=\sum_{p,q}f_{pq}$ as in~\eqref{UD}. The~following proposition shows that the right-hand side always makes sense.

\begin{proposition} \label{PD}
The series in \eqref{UG} converges absolutely and uniformly on compact subsets of~$\Bn$,
for~any $t\in\CC$ and any $f$ \Mh on~$\Bn$.
\end{proposition}

\begin{proof} For each nonnegative integer~$m$, the right-hand side of \eqref{UG} with $t=m$
is just the Peter-Weyl decomposition \eqref{UD} of the \Mh function $(I+\square)^m f$;
consequently, $\sum_{p,q}[4pq+(2n-2)(p+q)+1]^t |f_{pq}(\zeta)|$ converges uniformly on
compact subsets of~$\Bn$ for $t=m$, hence, a~fortiori, for any $t\in\CC$ with $\Re t\le m$.
Since $m\in\NN$ was arbitrary, the claim follows.
\end{proof}

We~proceed by recalling the formulas for the norms $\|f\|_s$, $s>-1$, in~terms of the Peter-Weyl
decomposition~\eqref{UD}. Employing again the polar coordinates $z=r\zeta$ on~$\CC^n$, from
$$ dz = \frac{2\pi^n}{\Gamma(n)} r^{2n-1} \,dr \,d\sigma(\zeta)  $$
we~obtain, for $f$ as in~\eqref{UD} and $s>-1$, using the previous notation $f_r(\zeta):=f(r\zeta)$,
\begin{align}
\|f\|^2_s &= \frac{(s+1)_n}{\pi^n} \frac{2\pi^n}{\Gamma(n)} \int_0^1 \int_\pBn
 |f(r\zeta)|^2 (1-r^2)^s r^{2n-1} \,d\sigma(\zeta) \,dr \nonumber \\
&= \frac{2(s+1)_n}{\Gamma(n)} \int_0^1 (1-r^2)^s r^{2n-1} \|f_r\|^2_\pBn \, dr \nonumber \\
&= \frac{2(s+1)_n}{\Gamma(n)} \int_0^1 (1-r^2)^s r^{2n-1} \sum_{p,q} \|(f_{pq})_r\|^2_\pBn \, dr \quad\text{by \eqref{UA}} \nonumber \\
&= \sum_{p,q} \frac{2(s+1)_n}{\Gamma(n)} \int_0^1 (1-r^2)^s r^{2n-1} \|S_{pq}(r^2) r^{p+q} f_{pq}\|^2_\pBn \, dr \quad\text{by \eqref{UC}} \label{UZ} \\
&= \sum_{p,q} \frac{(s+1)_n}{\Gamma(n)} \|f_{pq}\|^2_\pBn \int_0^1 (1-t)^s t^{n-1+p+q} S_{pq}(t)^2 \,dt \nonumber \\
&= \sum_{p,q} c_{pq}(s) \|f_{pq}\|^2_\pBn \quad\text{with $c_{pq}(s)$ as in \eqref{TI}}, \nonumber
\end{align}
which is the formula \eqref{TH} from the Introduction. Here the interchange of the integration
and summation is justified by the nonnegativity of the integrand.

One~more fact about \Mh functions which we will need is their \emph{invariant mean-value property}.
Recall that for each $z\in\Bn$, $z\neq0$, there is a unique $\phi_z\in\autBn$ which interchanges
$z$ and the origin~0; explicitly,
$$ \phi_z(w) = \frac{z-P_z w - \sqrt{1-|z|^2}(w-P_z w)}{1-\spr{w,z}}, \qquad P_z w:=\frac{\spr{w,z}}{|z|^2}z.  $$
For $z=0$, we~set $\phi_0(w):=-w$. For any $f$ \Mh on~$\Bn$, $z\in\Bn$ and $0<r<1$, one~then~has
\[ f(z) = \int_\pBn f(\phi_z(r\zeta)) \, d\sigma(\zeta) .  \label{UM}  \]
In~other words, $f(z)$ is equal to its mean value over any Moebius sphere centered at~$z$.

\section{The coefficients $c_{pq}(s)$} \label{Sec3}
The~following theorem is the crux of the proof of our main result.

\begin{theorem} \label{PE}
For $p,q\ge1$ and $s>-n-1$,
\[ \begin{aligned}
c_{pq}(s) &= \frac{\Gamma(s+n+1)\Gamma(s+1+2n)}{\Gamma(n)^3\Gamma(2s+2n+2)} (p)_n (q)_n \\
&\hskip2em
 \int_0^1 \int_0^1 x^{p-1} y^{q-1} (1-x)^{n+s} (1-y)^{n+s} \FF21{s+1,n+s+1}{2n+2s+2}{1-xy} \,dx\,dy.
\end{aligned} \label{VB} \]
\end{theorem}

We~actually present two proofs: the first one relies on a direct manipulation of the ${}_2\!F_1$ functions,
while the second one, inspired by the proof of Theorem~3.1 in Ureyen~\cite{Ur}, uses the trick of replacing
one of the $S_{pq}$ in~\eqref{TI} by its integral representation.

\begin{proof}[First proof of Theorem~\ref{PE}]
We~first need two lemmas. Let $\DD$ denote the unit disc in~$\CC$.

\begin{lemma}\label{hyplemma}
Let $\alpha>0$, $c>0$, $\sigma:=c-a-b+\alpha>0$, $\epsilon\in(0,1)$ and let $f$ be any function on
$\DD\setminus[-1,0]$ of the form
$$ f(z) = \sum_{k\in\NN} f_k z^k + \sum_{k\in\epsilon+\NN} f_k z^k =: g_1(z)+z^\epsilon g_2(z) $$
such that $g_1,g_2$ are holomorphic in~$\DD$. Then for any $x\in\DD\setminus[-1,0]$,
$$ \int_0^1 t^{c-1}(1-t)^{\alpha-1} \FF21{a,b}ct f(x(1-t)) \,dt = \Gamma(c) f\nadseb{\alpha,\sigma}{\sigma+a,\sigma+b}x,  $$
where
$$ f\nadseb{a,b}{c,d}x := \sum_{k\in\NN\cup(\epsilon+\NN)} f_k \frac{\Gamma(a+k)\Gamma(b+k)}{\Gamma(c+k)\Gamma(d+k)}x^k. $$
\end{lemma}

\begin{proof} Integrate term by term.
\end{proof}

\begin{lemma}\label{2F1/2F1} For $n\not\in \mathbb{Z}$ and for all $t\in(0,1)$:
\begin{align*}
\frac{\FF21{p,q}{p+q+n}t}{\FF21{p,q}{p+q+n}1} = \Big(\sum_{k\in\NN}-\sum_{k\in n+\NN}\Big)
 \frac{\Gamma(p+k)\Gamma(q+k)\Gamma(1-n)}{\Gamma(p)\Gamma(q)\Gamma(1-n+k)\Gamma(1+k)} (1-t)^k.
\end{align*}
\end{lemma}

\begin{proof} Immediate from the formula for analytic continuation of the ${}_2\!F_1$
function around $z=1$, cf.~\cite[formula~(1) in~\S2.10]{BE}.
\end{proof}

It~is now convenient to temporarily work with $n\notin\ZZ$, and in the end we take the limit~$n\to\NN$.
We~use the definition $(\nu)_n=\Gamma(\nu+n)/\Gamma(\nu)$ for the Pochhammer symbol $(\nu)_n$ with $\nu>0$ and arbitrary $n\ge0$.

\begin{corollary}\label{Cspqforcor} Let $n\notin\ZZ$, $n>0$. Then
\begin{align} \label{Cspqfor}
& c_{pq}(s) = \frac{(s+1)_n}{\Gamma(n)} \Big(\sum_{k\in\NN}-\sum_{k\in n+\NN}\Big) \\
& \frac{\Gamma(p+k)\Gamma(q+k)\Gamma(p+n)\Gamma(q+n)\Gamma(1-n)\Gamma(s+k+1)\Gamma(n+s+k+1)}{\Gamma(n)\Gamma(1-n+k)\Gamma(1+k)\Gamma(p)\Gamma(q)\Gamma(p+n+s+1+k)\Gamma(q+n+s+1+k)}.
 \nonumber
\end{align}
\end{corollary}

\begin{proof} Remember that
$$ c_{pq}(s) = \frac{(s+1)_n}{\Gamma(n)} \int_0^1 t^{p+q+n-1}(1-t)^s \left(\frac{\FF21{p,q}{p+q+n}t}{\FF21{p,q}{p+q+n}1}\right)^2 \,dt. $$
Denote
$$ f(x):=\Big(\sum_{k\in\NN}-\sum_{k\in n+\NN}\Big) \frac{\Gamma(p+k)\Gamma(q+k)\Gamma(1-n)}{\Gamma(p)\Gamma(q)\Gamma(1-n+k)\Gamma(1+k)} x^k. $$
Using Lemma \ref{2F1/2F1} we have
$$ c_{pq}(s)=\frac{(s+1)_n}{\Gamma(n)} \int_0^1 t^{p+q+n-1} (1-t)^s \frac{\FF21{p,q}{p+q+n}t}{\FF21{p,q}{p+q+n}1} f(1-t)\,dt. $$
For a moment we insert a variable $x$ into the picture and evaluate the ${}_2\!F_1(1)$ in Gamma functions. Observe that
$$ c_{pq}(s) = \lim_{x\nearrow 1} \frac{(s+1)_n}{\Gamma(n)}
 \frac{\Gamma(p+n)\Gamma(q+n)}{\Gamma(p+q+n)\Gamma(n)}
 \int_0^1 t^{p+q+n-1}(1-t)^s \FF21{p,q}{p+q+n}t f((1-t)x) \,dt , $$
the interchange of the limit and the integral being justified by the Lebesgue Dominated Convergence Theorem,
since for $n>0$ the ${}_2\!F_1$ function is bounded on~$[0,1]$.
By~Lemma \ref{hyplemma} we thus get
\begin{align*}
c_{pq}(s) &= \lim_{x\nearrow 1} \frac{(s+1)_n}{\Gamma(n)}
 \frac{\Gamma(p+n)\Gamma(q+n)}{\Gamma(n)} f\nadseb{s+1,n+s+1}{p+n+s+1,q+n+s+1}x \\
&= \frac{(s+1)_n}{\Gamma(n)}
 \frac{\Gamma(p+n)\Gamma(q+n)}{\Gamma(n)} f\nadseb{s+1,n+s+1}{p+n+s+1,q+n+s+1}1.
\end{align*}
This is exactly (\ref{Cspqfor}).
\end{proof}

We~are now ready to finish the first proof of Theorem~\ref{PE}.
Assume first that $n\notin\ZZ$, $n>0$; in the end we take the limit $n\to\NN$.
Using the integral representations:
\begin{align*}
\frac{\Gamma(p+k)\Gamma(n+s+1)}{\Gamma(p+n+s+1+k)}&=\int_0^1 x^{p+k-1}(1-x)^{n+s} \,dx,\\
\frac{\Gamma(q+k)\Gamma(n+s+1)}{\Gamma(q+n+s+1+k)}&=\int_0^1 y^{q+k-1}(1-y)^{n+s} \,dy,
\end{align*}
in~formula \eqref{Cspqfor} we~get
$$ c_{pq}(s) = \frac{(s+1)_n}{\Gamma(n)} \int_0^1 \int_0^1 x^{p-1} y^{q-1} (1-x)^{n+s} (1-y)^{n+s} g(xy)\,dx\,dy $$
where
\begin{align*}
g(z)&:=\Big(\sum_{k\in\NN}-\sum_{k\in n+\NN}\Big) \frac{(p)_n(q)_n\Gamma(1-n)\Gamma(s+k+1)\Gamma(n+s+k+1)}
 {\Gamma(n)\Gamma(1-n+k)\Gamma(1+k)\Gamma(n+s+1)^2} z^k \\
&=(p)_n(q)_n \frac{\Gamma(s+1)\Gamma(s+1+2n)}{\Gamma(n)^2\Gamma(2s+2n+2)} \FF21{s+1,n+s+1}{2n+2s+2}{1-z},
\end{align*}
where the last equality was obtain by applying Lemma \ref{2F1/2F1} with $p$, $q$, and $t$ replaced by $s+1$, $n+s+1$, and $1-z$, respectively.

The resulting formula makes sense also for $n\in \mathbb{N}$.
Thus now we take our limit, which yields~\eqref{VB}.
\end{proof}

\begin{proof}[Second proof of Theorem~\ref{PE}]
The~following alternative proof was inspired by the paper~\cite{Ur}.
Remember that
\[ \label{E69} c_{pq}(s) = \frac{(s+1)_n}{\Gamma(n)} \int_0^1 t^{p+q+n-1} (1-t)^s \left(\frac{\FF21{p,q}{p+q+n}t}{\FF21{p,q}{p+q+n}1}\right)^2 \,dt. \]
Representing one of the ${}_2\!F_1$ functions~as
$$ \frac{\FF21{p,q}{p+q+n}t}{\FF21{p,q}{p+q+n}1} = \frac{(p)_n}{\Gamma(n)} \int_0^1 r^{p-1}(1-r)^{q+n-1}(1-rt)^{-q} \,dr $$
and inserting this into (\ref{E69}) we~get
\begin{align}
c_{pq}(s) &= \frac{(s+1)_n}{\Gamma(n)} \frac{(p)_n\Gamma(p+n)\Gamma(q+n)}{\Gamma(n)^2\Gamma(p+q+n)} \nonumber \\
& \hskip4em \int_0^1 \int_0^1 t^{p+q+n-1}(1-t)^s r^{p-1}(1-r)^{q+n-1}(1-rt)^{-q} \FF21{p,q}{p+q+n}t \,dr \,dt \nonumber \\
\label{E71} &= \frac{(s+1)_n}{\Gamma(n)} \frac{(p)_n\Gamma(p+n)\Gamma(q+n)}{\Gamma(n)^2\Gamma(p+q+n)} \\
& \hskip4em \int_0^1 r^{p-1}(1-r)^{q+n-1} \int_0^1 t^{p+q+n-1}(1-t)^s (1-rt)^{-q} \FF21{p,q}{p+q+n}t \,dt \,dr. \nonumber
\end{align}
On~the inner integral we use Lemma~2.10 from \cite{Ur} to get
\begin{multline*}
\int_0^1 t^{p+q+n-1} (1-t)^s (1-rt)^{-q} \FF21{p,q}{p+q+n}t \,dt =
 \GG{p+q+n,s+1,n+s+1}{p+n+s+1,q+n+s+1} \\
\times (1-r)^{q}\FF32{s+1,q,n+s+1}{p+n+s+1,q+n+s+1}{\frac{r}{r-1}} \\
=\GG{p+q+n}{q,p+n} (1-r)^{-q} \int_0^1 \int_0^1 (1-u)^s u^{p+n-1} y^{q-1} (1-y)^{n+s} \Big(1-\frac{(1-u)yr}{r-1}\Big)^{-(n+s+1)} \,du \,dy .
\end{multline*}
Inserting this into (\ref{E71}) and rearranging we obtain
$$ c_{pq}(s) = \frac{(s+1)_n}{\Gamma(n)} \frac{(p)_n(q)_n}{\Gamma(n)^2}
 \int_0^1 \int_0^1 \int_0^1 \frac{r^{p-1}(1-r)^{2n+s}(1-u)^{s} u^{p+n-1}y^{q-1}(1-y)^{n+s}}{(1-r(1-(1-u)y))^{n+s+1}} \,du \,dr \,dy. $$
We~now change the integration variable $u$ to $u=x/r$ so that $du=dx/r$. Then we interchange the order of integration
$$ \int_0^1 \int_0^r \,dx \,dr \quad\to\quad  \int_0^1 \int_x^1 \,dr \,dx, $$
and also change the variable $r$ to $r\to (1-x)r+x$. Overall we obtain
\begin{align*}
c_{pq}(s) &= \frac{(s+1)_n}{\Gamma(n)} \frac{(p)_n(q)_n}{\Gamma(n)^2} \int_0^1 \int_0^1 \int_0^1
 \frac{x^{p+n-1}y^{q-1}(1-x)^{n+s}(1-y)^{n+s} r^s(1-r)^{s+2n}}{(x+r(1-x))^{-(n+s+1)}(1-r(1-v))^{-(n+s+1)}} \,dr \,dx \,dy \\
&=\frac{(s+1)_n}{\Gamma(n)} \frac{(p)_n(q)_n \Gamma(s+1)\Gamma(s+2n+1)}{\Gamma(n)^2\Gamma(2s+2n+2)}
 \int_0^1 \int_0^1 x^{p+n-1}y^{q-1}(1-x)^{n+s}(1-y)^{n+s} \\
& \hskip4em x^{-n-s-1} \FE1{s+1;n+s+1,n+s+1}{2s+2n+2}{1-\frac1x,1-y} \,dx \,dy.
\end{align*}
A~well known transform of the Appell $F_1$ function
$$ \FE1{a;b_1,b_2}{b_1+b_2}{x,y} = (1-x)^{-a} \FF21{a,b_2}{b_1+b_2}{\frac{x-y}{x-1}} $$
(see \cite[formula~(1) in~\S5.10]{BE}) implies
$$ x^{-n-s-1} \FE1{s+1;n+s+1,n+s+1}{2s+2n+2}{1-\frac1x,1-y} = x^{-n}\FF21{s+1,n+s+1}{2s+2n+2}{1-xy}. $$
Overall, we~get
\begin{multline*}
c_{pq}(s) = \frac{(s+1)_n}{\Gamma(n)} \frac{(p)_n(q)_n \Gamma(s+1)\Gamma(s+2n+1)}{\Gamma(n)^2\Gamma(2s+2n+2)} \\
 \int_0^1 \int_0^1 x^{p-1} y^{q-1} (1-x)^{n+s} (1-y)^{n+s} \FF21{s+1,n+s+1}{2s+2n+2}{1-x y} \,dx \,dy,
\end{multline*}
as claimed.
\end{proof}

We~remark that for $s>-1$, using the identities
$$ (p)_n x^{p-1} = \partial_x^n x^{p+n-1}, \qquad (q)_n y^{q-1} = \partial_y^n y^{q+n-1}, $$
and integrating by parts $n$ times in both integrals, one~can rewite~\eqref{VB}~as
$$ c_{pq}(s) = \frac{\Gamma(s+1+2n)\Gamma(s+n+1)}{\Gamma(n)^3\Gamma(2s+2n+2)}
 \int_0^1 \int_0^1 x^{p+n-1} y^{q+n-1} f(x,y) \,dx \,dy, $$
where
$$ f(x,y) = \frac{\partial^{2n}}{\partial x^n\partial y^n} (1-x)^{n+s}(1-y)^{n+s} \FF21{s+1,n+s+1}{2n+2s+2}{1-xy}. $$
Though the function $f(x,y)$ here is more complicated than the single ${}_2\!F_1$ in~\eqref{VB},
this formula has the advantage of being valid also for $p=0$ or $q=0$,
as~is readily checked by direct verification.

\begin{corollary} \label{PC}
For each $s>-n-1$, there exist constants $0<c_s<C_s<+\infty$ such that
\[ c_s \le (p+1)^{s+1}(q+1)^{s+1} c_{pq}(s) \le C_s \qquad\forall p,q\ge0.  \label{VG} \]
\end{corollary}

\begin{proof} For $p=0$, so~that $S_{pq}\equiv1$, direct evaluation of~\eqref{TI} yields
\[ c_{0q}(s) = \frac{\Gamma(q+n)\Gamma(n+s+1)}{\Gamma(n+s+q+1)\Gamma(n)}  \label{VC}  \]
which is $\sim q^{-s-1}$ as $q\to+\infty$, by Stirling's formula; thus \eqref{VG} holds
for $p=0$. Similarly for $q=0$.

For $p,q\ge1$, the ${}_2\!F_1$ function in~\eqref{VB} is positive and continuous on the closed
interval~$[0,1]$, with the value 1 at 0 and the finite value $\frac{\Gamma(2n+2s+2)\Gamma(n)}
{\Gamma(2n+s+1)\Gamma(n+s+1)}$ at~1. This is clear for $s>-1$, while for general $s>-n-1$,
the~positivity follows from Euler's formula~\cite[formula~(23) in~\S2.1]{BE}
$$ \FF21{s+1,n+s+1}{2n+2s+2}t = (1-t)^n \FF21{2n+s+1,n+s+1}{2n+2s+2}t . $$
Thus the value of the expression
$$ (p)_n (q)_n \int_0^1 \int_0^1 x^{p-1} y^{q-1} (1-x)^{n+s} (1-y)^{n+s} \FF21{s+1,n+s+1}{2n+2s+2}{1-xy} \,dx\,dy $$
in~\eqref{VB} is bounded from above as well as from below by constant multiples~of
\begin{multline*}
(p)_n (q)_n \int_0^1 \int_0^1 x^{p-1} y^{q-1} (1-x)^{n+s} (1-y)^{n+s} \,dx\,dy \\
= (p)_n (q)_n \GG{p,n+s+1}{n+s+p+1} \GG{q,n+s+1}{n+s+q+1} \\
= \GG{n+p,n+s+1}{n+s+p+1} \GG{n+q,n+s+1}{n+s+q+1} .
\end{multline*}
Again, the first factor is $\asymp(p+1)^{-s-1}$ and the second factor is $\asymp(q+1)^{-s-1}$ by Stirling's formula.
\end{proof}

\section{\Mh Besov spaces} \label{Sec4}
We~are now ready to prove our first main result.
With the estimate \eqref{TG} in hands, the theorem below follows much in the same way as for Theorem~9 in~\cite{EY2}.

\begin{theorem} \label{PA}
For $f$ \Mh on $\Bn$, $n\ge2$, and $s>-n-1$, the following assertions are equivalent.
\begin{itemize}
\item[(a)] $f\in\MM_s$;
\item[(b)] $\sum_{j_1,\dots,j_m=1}^{2n^2} \|X_{j_1}\dots X_{j_m}f\|^2_{s+m} + |f(0)|^2 <+\infty$
for some (equivalently, any) $m=0,1,2,\dots$, $m>-s-1$;
\item[(c)] $\|(I+\square)^{t/2}f\|^2_{s+t}<+\infty$ for some (equivalently, any) $t>-s-1$;
\item[(d)] $\|(I+\square)^{-(s+1)/2}f\|^2_\Hardy<+\infty$.
\end{itemize}
Furthermore, the quantities above are all equivalent to $\|f\|^2_s$.
\end{theorem}

\begin{proof} We~will show that for $f=\sum_{p,q}f_{pq}$, $f_{pq}\in\bhpq$, \Mh on~$\Bn$,
$\|f\|_s^2$ as well as the expressions in the items (b)--(d) are all equivalent~to
$$ \sum_{p,q} \frac{\|f_{pq}\|^2_\pBn}{(p+1)^{s+1}(q+1)^{s+1}} ;  $$
this will settle the claim.

(a) By~\eqref{TH} and its analytic continuation to $s>-n-1$ \cite{EY2}, we~have
$$ \|f\|_s^2 = \sum_{p,q} c_{pq}(s) \|f_{pq}\|^2_\pBn  $$
(for $-n-1<s\le-1$, this is actually the definition of~$\|f\|_s$).
The~assertion is thus immediate from~\eqref{TG}.

(b) Since $-\oL_{jk}$ is the adjoint of $L_{jk}$ in $L^2(\pBn,d\sigma)$,
we~have for any $g\in C^2(\pBn)$
$$ \sum_{j=1}^{2n^2} \|X_j g\|^2_\pBn = - \sum_{j,k=1}^n \spr{(L_{jk}\oL_{jk}+\oL_{jk} L_{jk})g,g}_\pBn
 = \spr{\square g,g}_\pBn . $$
If~$g=\sum_{p,q}g_{pq}$, $g_{pq}\in\chpq$, is~the Peter-Weyl decomposition~\eqref{UA} of~$g$,
we~thus have by~\eqref{UF}
$$ \sum_{j=1}^{2n^2} \|X_j g\|^2_\pBn = \sum_{p,q} [4pq+(2n-2)(p+q)] \|g_{pq}\|^2_\pBn.  $$
Iterating this formula, we~obtain
\[ \sum_{j_1,j_2,\dots,j_m=1}^{2n^2} \|X_{j_1}X_{j_2}\dots X_{j_m}g\|^2_\pBn
 = \sum_{p,q} [4pq+(2n-2)(p+q)]^m \|g_{pq}\|^2_\pBn , \]
for any $m=0,1,2,\dots$.
Applying this to $g=f_r$, where again $f_r(\zeta):=f(r\zeta)$, $0<r<1$, yields,
thanks to~\eqref{UC},
$$ \sum_{j_1,j_2,\dots,j_m=1}^{2n^2} \|X_{j_1}X_{j_2}\dots X_{j_m}f_r\|^2_\pBn
 = \sum_{p,q} [4pq+(2n-2)(p+q)]^m r^{2(p+q)} S^{pq}(r^2)^2 \|f_{pq}\|^2_\pBn . $$
As~$X_j$, being tangential, do~not act on the $r$ variable, we~also have
$$ X_{j_1}X_{j_2}\dots X_{j_m}f_r = (X_{j_1}X_{j_2}\dots X_{j_m}f)_r . $$
Hence for any $s>-1$,
\begin{align*}
\|X_{j_1}X_{j_2}\dots X_{j_m}f\|^2_s
&= \frac{(s+1)_n}{\pi^n} \int_0^1 \frac{2\pi^n}{\Gamma(n)} \|(X_{j_1}X_{j_2}\dots X_{j_m}f)_r\|^2_\pBn (1-r^2)^s r^{2n-1}\,dr \\
&= \frac{(s+1)_n}{\Gamma(n)} \sum_{p,q} [4pq+(2n-2)(p+q)]^m \|f_{pq}\|^2_\pBn \int_0^1 S^{pq}(t)^2 t^{n-1} (1-t)^s \,dt \\
&= \sum_{p,q} [4pq+(2n-2)(p+q)]^m \|f_{pq}\|^2_\pBn c_{pq}(s).
\end{align*}
Denoting temporarily
$$ d_{pq} := \begin{cases} 1 & \text{if } p=q=0 \\
[4pq+(2n-2)(p+q)]^m c_{pq}(s+m) \quad & \text{if } p+q>0 \end{cases} $$
we~thus have
$$ \sum_{j_1,\dots,j_m=1}^{2n^2} \|X_{j_1}\dots X_{j_m}f\|^2_{s+m} + |f(0)|^2
 = \sum_{p,q} d_{pq} \|f_{pq}\|^2_\pBn . $$
However, for $p+q>0$ and $n\ge2$ clearly
$$ [4pq+(2n-2)(p+q)] \asymp (p+1)(q+1).  $$
Hence by \eqref{TG} $d_{pq}\asymp(p+1)^{-s-1}(q+1)^{-s-1}$ for all $p,q\ge0$,
and the assertion follows.

(c) By~\eqref{UG} and~\eqref{UA}, we~have for any $0<r<1$,
$$ \|(I+\square)^{t/2}f_r\|^2_\pBn = \sum_{p,q} [4pq+(2n-2)(p+q)+1]^t r^{2(p+q)} S_{pq}(r^2)^2 \|f_{pq}\|^2_\pBn . $$
The same calculation as in part (b) therefore shows that for any $s>-1$,
$$ \|(I+\square)^{t/2}f\|^2_s = \sum_{p,q} [4pq+(2n-2)(p+q)+1]^t c_{pq}(s) \|f_{pq}\|^2_\pBn . $$
Since $[4pq+(2n-2)(p+q)+1]^t c_{pq}(s+t)\asymp(p+1)^{-s-1}(q+1)^{-s-1}$ by~\eqref{TG},
the assertion follows.

(d) Again by~\eqref{UG} and~\eqref{UA}, for any $0<r<1$,
$$ \|(I+\square)^{-(s+1)/2}f_r\|^2_\pBn = \sum_{p,q} [4pq+(2n-2)(p+q)+1]^{-s-1} r^{2(p+q)} S_{pq}(r^2)^2 \|f_{pq}\|^2_\pBn . $$
Letting $r\nearrow1$, we~have $r^{p+q}S_{pq}(r^2)\nearrow1$, so~by the Lebesgue Monotone Convergence Theorem
$$ \|(I+\square)^{-(s+1)/2}f_r\|^2_\Hardy = \sum_{p,q} [4pq+(2n-2)(p+q)+1]^{-s-1} \|f_{pq}\|^2_\pBn . $$
Since $[4pq+(2n-2)(p+q)+1]\asymp(p+1)(q+1)$, we~are done.
\end{proof}

Denoting, as~in the Introduction, for any real~$s$,
\[ \begin{aligned}
\wt\MM_s :&= \{ f\text{ \Mh on }\Bn: \;\|f\|_{\tilde s}:=\|(I+\square)^{-(s+1)/2}f\|^2_\Hardy<+\infty \} , \\
\MM_{\#,s} :&= \Big\{f=\sum_{p,q} f_{pq},\;f_{pq}\in\bhpq: \; \|f\|^2_{\#,s}:=\sum_{p,q} \frac{\|f_{pq}\|^2_\pBn}{(p+1)^{s+1}(q+1)^{s+1}} < +\infty \Big\},
\end{aligned} \label{VD} \]
we~have thus shown that
\[ \wt\MM_s=\MM_{\#,s} \quad\forall s\in\RR, \qquad \wt\MM_s=\MM_{\#,s}=\MM_s \quad \forall s>-n-1,  \label{VE} \]
with equivalent norms.

\section{\Mh Sobolev spaces} \label{Sec5}
In~addition to our operators $X_j$, $j=1,\dots,2n^2$, from \eqref{TD} in the Introduction,
and the complex normal derivative \eqref{UR} which we now denote by~$X_0$:
$$ X_0 := \cR , $$
we~will now need also the real normal derivative
$$ \cN = \sum_{j=1}^n \Big(z_j\frac\partial{\partial z_j} + \oz_j\frac\partial{\partial\oz_j}\Big), $$
or, in~the polar coordinates $z=r\zeta$, $\cN=r\partial/{\partial r}$.

Introduce the quantities
\[ c_{pq,k}(s) := \int_0^1 t^{n-1} (1-t)^s [(2t\partial/\partial t)^k (t^{\frac{p+q}2} S_{pq}(t))]^2 \, dt , \label{WA} \]
so~that
\[ \frac{(s+1)_n}{\Gamma(n)} c_{pq,0}(s) = c_{pq}(s). \label{ccc}  \]
Clearly
\[ \begin{aligned}
c_{00,k}(s) &= 0 \qquad\text{for }k>0 , \\
c_{0p,k}(s) &= c_{p0,k}(s) = p^{2k} \frac{\Gamma(n+p)}{(s+1)_{n+p}} \asymp (p+1)^{2k-s-1} \qquad\text{for }p>0 ,
\end{aligned} \label{WB} \]
by Stirling's formula.
The~significance of the quantities $c_{pq,k}(s)$ stems from the following proposition.

\begin{proposition} \label{PF}
For any $m=0,1,2,\dots,n$, $s>-1$ and $f$ \Mh on~$\Bn$, $n\ge2$, the following assertions are equivalent.
\begin{itemize}
\item[(a)] $\sum_{l=0}^m \sum_{k=0}^l \sum_{j_1,\dots,j_{l-k}=0}^{2n^2} \|\cN^k X_{j_1}\dots X_{j_{l-k}}f\|^2_s <+\infty$;
\item[(b)] $\sum_{l=0}^m \sum_{k=0}^l (p+q)^{2(l-k)} c_{pq,k}(s) \|f_{pq}\|^2_\pBn <+\infty$;
\item[(c)] $f\in W^m(\Bn,d\mu_s)$, the weighted Sobolev space of order $m$ on $\Bn$ with respect to~$d\mu_s$.
\end{itemize}
Furthermore, the quantities in (a), (b) are equivalent to the squared norm of $f$ in~$W^m(\Bn,d\mu_s)$.
\end{proposition}

Here (b) means that if $c_{pq,k}(s)=+\infty$ for some $p,q,k$, then $f_{pq}\equiv0$,
and $c_{pq,k}(s)\|f_{pq}\|^2_\pBn$ is~to be interpreted as 0 for such $p,q,k$.
Note also that the tangential operators $X_j$ automatically commute with the normal
derivative~$\cN$, so~the order of operators in (a) is irrelevant.

\begin{proof} (a)$\iff$(b) Note that for $f=\sum_{p,q}f_{pq}$, $f_{pq}\in\bhpq$, as~in~\eqref{UD},
we~have $\cN f=\sum_{p,q}\cN f_{pq}$, because the operator $\cN$ acts only on the $r$ variable in
$f_{pq}(r\zeta)=r^{p+q}S_{pq}(r^2)f_{pq}(\zeta)$, thus preserving the Peter-Weyl components with
respect to the $\zeta$ variable. By~our calculations in the proof of Theorem~\ref{PA} and as in~\eqref{UZ},
we~thus have
\begin{align*}
\sum_{j_1,\dots,j_{l-k}=0}^{2n^2} \|X_{j_1}\dots X_{j_{l-k}}\cN^k f\|^2_s
 &= \sum_{p,q} [(p+q)(p+q+2n-2)]^{l-k} \|\cN^k f_{pq}\|^2_s \\
&= \frac{(s+1)_n}{\Gamma(n)} \sum_{p,q} [(p+q)(p+q+2n-2)]^{l-k} c_{pq,k}(s) \|f_{pq}\|^2_s ,
\end{align*}
thus proving the equivalence of (a) and~(b).

(c)$\implies$(a) By~the Leibniz rule,
$$ X_{j_1}\dots X_{j_{l-k}}\cN^k f(z)=\sum_{|\alpha|+|\beta|\le l}
P_{j_1,\dots,j_{l-k};k;\alpha\beta}(z) \partial^\alpha\dbar{}^\beta f(z) $$
with some coefficient functions
$P_{j_1,\dots,j_{l-k};k;\alpha\beta}$ which are bounded on~$\Bn$ (in~fact --- they are polynomials).
The~claim is thus immediate from the triangle inequality.

(a)$\implies$(c) Observe that at any $z\in\Bn$, $z\neq0$, the tangential vector fields~$X_j$,
$j=0,\dots,2n^2$, span (very redundantly) the entire real tangent space to the sphere $|z|\pBn$;
thus together with~$\cN$, they span the whole tangent space at $z$ in~$\Bn$. It~follows that,
for~any $l\in\NN$, the derivatives $\partial^\alpha\dbar{}^\beta f$, $|\alpha|+|\beta|=l$, can be
expressed as linear combinations of the derivatives $X_{j_1}\dots X_{j_{l-k}}\cN^k f$, $0\le k\le l$,
$0\le j_1,\dots,j_{l-k}\le2n^2$; furthermore, the~coefficients of these linear combinations can be
chosen to be bounded outside any neighborhood of the origin $z=0$ (where all the $X_j$ as well as
$\cN$ vanish). In~other words, if~we momentarily denote by~$\chi$ the characteristic function of
the annular region $\frac12<|z|<1$, then for any $s>-1$ and $l\in\NN$,
$$ \sum_{|\alpha|+|\beta|\le l} \|\chi\partial^\alpha\dbar{}^\beta f\|^2_s \lesssim
 \sum_{k=0}^l \sum_{j_1,\dots,j_{l-k}=0}^{2n^2} \|X_{j_1}\dots X_{j_{l-k}}\cN^k f\|^2_s,  $$
and, hence,
\[ \sum_{l=0}^m \sum_{|\alpha|+|\beta|\le l} \|\chi\partial^\alpha\dbar{}^\beta f\|^2_s \lesssim
 \sum_{l=0}^m \sum_{k=0}^l \sum_{j_1,\dots,j_{l-k}=0}^{2n^2} \|X_{j_1}\dots X_{j_{l-k}}\cN^k f\|^2_s .  \label{WC}  \]
To~treat the remaining region $|z|<\frac12$, we~use a ``subharmonicity'' argument.
Let~$\theta$ be any smooth function on~$\Bn$ supported on~$|z|<\frac12$, depending only on~$|z|$
and of total mass~one. By~\eqref{UM}, for~any $a\in\Bn$ and $f$ \Mh on~$\Bn$,
$$ f(a) = \int_\Bn f(\phi_a(z)) \, \theta(z) \, dz.  $$
Changing the variable from $z$ to $\phi_a(z)$ gives
$$ f(a) = \int_\Bn f(z) \theta(\phi_a(z)) J(a,z) \, d\mu_s(z),  $$
with
$$ J(a,z) := |\operatorname{Jac}_{\phi_a}(z)|^2 \frac{\pi^n}{(s+1)_n} (1-|z|^2)^s  $$
smooth on $\Bn\times\Bn$. Hence, for any multi-indexes $\alpha$ and~$\beta$,
$$ \partial^\alpha \dbar{}^\beta f(a) = \int_\Bn f(z)
 \frac{\partial^{|\alpha|+|\beta|}}{\partial a^\alpha\partial\overline a{}^\beta}
 [\theta(\phi_a(z) J(a,z)] \, d\mu_s(z)  . $$
Since $\theta$ is supported on~$|z|<\frac12$, the integration is in fact only over
$\{z:|\phi_a(z)|<\frac12\}$; for~$|a|<\frac12$, this implies that $|z|<\rho$ for
some $\rho<1$. Since $J(a,z)$ is bounded for $|a|\le\frac12$ and $|z|\le\rho$,
we~thus~get
$$ |\partial^\alpha \dbar{}^\beta f(a)| \le C_{\alpha\beta} \int_{|z|<\rho} |f(z)| \,d\mu_s(z) \le C_{\alpha\beta} \|f\|_s $$
with some finite $C_{\alpha\beta}$ independent of $|a|\le\frac12$ and~$f$. Hence
$$ \|(1-\chi)\partial^\alpha \dbar{}^\beta f\|_s \le C_{\alpha\beta} \|f\|_s \|1-\chi\|_s \le C_{\alpha\beta}\|f\|_s,  $$
and, consequently,
$$ \sum_{l=0}^m \sum_{|\alpha|+|\beta|\le l} \|(1-\chi)\partial^\alpha\dbar{}^\beta f\|^2_s \lesssim \|f\|^2_s . $$
Combining this with~\eqref{WC}, the claim follows.
\end{proof}

\begin{remark} The~above proof is modeled on the proof of the implication (e)$\implies$(f) of Theorem~14 in~\cite{BEY},
rectifying a small error we made there. (We~used for $\theta$ just the characteristic function of $|z|<\frac14$,
which is of course not smooth at $|z|=\frac14$.)  \qed
\end{remark}

We~now need the following analogue of Corollary~\ref{PC} for the quantities~$c_{pq,k}(s)$.

\begin{proposition} \label{PG}
For each $s>-1$, $n\ge2$ and $k=0,1,\dots,n$,
\[ \begin{aligned} c_{00,k}(s) &\asymp \begin{cases} 1 \quad & k=0,\\0 & k>0,\end{cases} \\
c_{pq,k}(s) &\asymp [(p+1)(q+1)]^{2k-s-1} \qquad\text{for all }p+q>0.
\end{aligned} \label{WD}  \]
\end{proposition}

For the proof, we~first need a lemma and a proposition.
Denote, for $p+q+n>0$, $k\ge0$ and $s>-1$,
$$ I_{pqs}(n,k) := \int_0^1 t^{p+q+n-1+k} (1-t)^s \FF21{p,q}{p+q+n}t ^2 \, dt.  $$
Thus, for all $p,q\ge0$ and $n\ge1$,
\[ c_{pq,0}(s) = \GG{n+p,n+q}{n,n+p+q} ^2 I_{pqs}(n,0),  \label{mA} . \]

\begin{lemma} \label{Ap2} For $p+q>0$, $n\ge1$ and $0\le k\le n$,
\[ c_{pq,k}(s) \asymp (p+q+1)^{2k} \GG{n+p,n+q}{n,n+p+q} ^2 I_{pqs}(n-k,k)  \label{mB}. \]
\end{lemma}

\begin{proof} For $pq=0$ and $k<n$, this is immediate from~\eqref{WB},
so~it is enough to consider $p+q\ge1$ and $0\le k\le n$. Recall that
\begin{align*}
c_{pq,k}(s) &:= \int_0^1 t^{n-1} (1-t)^s [(2t\partial/\partial t)^k (t^{\frac{p+q}2}S_{pq}(t))]^2 \,dt \\
&= \GG{n+p,n+q}{n,n+p+q}^2 \int_0^1 t^{n-1} (1-t)^s \Big[(2t\partial/\partial t)^k \Big(t^{\frac{p+q}2}\FF21{p,q}{p+q+n}t\Big)\Big]^2 \,dt , \\
c_{pq}(s) &:= \frac{(s+1)_n}{\Gamma(n)} c_{pq,0}(s).
\end{align*}
From the Taylor expansion for ${}_2\!F_1$, we~have
$$ (2t\partial/\partial t)^k \Big(t^{\frac{p+q}2}\FF21{p,q}{p+q+n}t\Big) = t^{\frac{p+q}2} \sum_{j=0}^\infty \frac{(p)_j(q)_j}{j!(n+p+q)_j} (p+q+2j)^k t^j . $$
Now
$$ \GG{p+q+n-k+j}{p+q+n+j} (p+q+2j)^k = \prod_{\ell=0}^{k-1} \frac{p+q+2j}{p+q+n-k+j+\ell} . $$
The $\ell$-th factor is a monotone function of $p+q\ge1$ (decreasing for $j>n-k+\ell$, increasing for $j<n-k+\ell$),
with limit 1 for $p+q\to+\infty$; its value thus always lies between 1 and $\frac{2j+1}{n-k+j+\ell+1}$.
The latter expression is again a monotone function of~$j\ge0$, with values between 2 and $\frac1{n-k+\ell+1}$.
Consequently, the whole product must lie between $2^k$ and $1/(n-k+1)_k$, and
$$ \frac1{(n-k+1)_k} \frac1{\Gamma(p+q+n-k+j)}\le \frac{(p+q+2j)^k}{\Gamma(p+q+n+j)} \le 2^k \frac1{\Gamma(p+q+n-k+j)} . $$
Multiplying by $\Gamma(p+q+n)$ yields
$$ \frac{(p+q+n-k)_k}{(n-k+1)_k} \frac1{(p+q+n-k)_j}\le \frac{(p+q+2j)^k}{(p+q+n)_j} \le 2^k (p+q+n-k)_k \frac1{(p+q+n-k)_j} . $$
Thus
\begin{multline*}
 t^{\frac{p+q}2} \frac{(p+q+n-k)_k}{(n-k+1)_k} \FF21{p,q}{p+q+n-k}t \\
 \le (2t\partial/\partial t)^k \Big(t^{\frac{p+q}2}\FF21{p,q}{p+q+n}t\Big) \\
 \le t^{\frac{p+q}2} 2^k (p+q+n-k)_k \FF21{p,q}{p+q+n-k}t,
\end{multline*}
and
\begin{multline*}
 \frac{(p+q+n-k)_k^2}{(n-k+1)_k^2} I_{pqs}(n-k,k) \\
 \le \GG{n,n+p+q}{n+p,n+q}^2 c_{pq,k}(s) \\
 \le 4^k (p+q+n-k)_k^2 I_{pqs}(n-k,k).
\end{multline*}
As $(p+q+n-k)_k\asymp(p+q+1)^k$, \eqref{mB} follows.
\end{proof}

\begin{proposition} \label{cc2}
For any $k\ge0$, $n\ge0$ and $s>-1$,
$$ I_{pqs}(n,k) \asymp I_{pqs}(n,0) \qquad\text{uniformly } \forall p+q\ge1.  $$
\end{proposition}

\begin{proof} Let
$$ \FF21{p,q}{p+q+n}t ^2 = \sum_{j=0}^\infty a_j t^j  $$
be the Taylor expansion; clearly $a_j\ge0$ $\forall j$. Integration term by term yields
$$ I_{pqs}(n,k) = \sum_{j=0}^\infty \GG{b+j+k,s+1}{b+j+k+s+1} a_j, \qquad b:=p+q+n. $$
Now
$$ \GG{b+j+k,s+1}{b+j+k+s+1} \Big/ \GG{b+j,s+1}{b+j+s+1} = \GG{b+j+k,b+j+s+1}{b+j,b+j+k+s+1} . $$
This is a positive function of $b+j\in[1,+\infty)$, with limit 1 at the infinity and positive value
$k!/(s+2)_k$ at $b+j=1$. Hence there exist constants $0<c_{ks}<C_{ks}<+\infty$ such that
$$ c_{ks} \le \GG{b+j+k,b+j+s+1}{b+j,b+j+k+s+1} \le C_{ks} \qquad \forall b\ge1, \forall j\ge0.  $$
Consequently,
$$ c_{ks} \le \frac{I_{pqs}(n,k)} {I_{pqs}(n,0)} \le C_{ks} \qquad \forall p+q\ge1, \forall n\ge0,  $$
proving the claim.
\end{proof}

\begin{corollary} For any $n\ge2$, $k=0,1,2,\dots,n$ and $s>-1$,
$$ c_{pq,k}(s) \asymp [(p+1)(q+1)]^{2k} c_{pq,0}(s), \qquad \text{uniformly for all }p+q>0.  $$
\end{corollary}

\begin{proof} For $k=0$ this is trivial, and for $k>0$, $pq=0$ the claim is immediate from the explicit calculation
$$ c_{p0,k}(s)=p^{2k} c_{p0,0}(s), \quad c_{0q,k}(s)=q^{2k}c_{0q,0}(s), $$
(cf.~\eqref{WB}), so~it is enough to consider $k>0$ and $p,q\ge1$. From \eqref{mB} and the last proposition,
\begin{align}
\frac{c_{pq,k}(s)}{(p+q+1)^{2k}} & \asymp \GG{n+p,n+q}{n,n+p+q} ^2 I_{pqs}(n-k,k) \nonumber \\
& \asymp \GG{n+p,n+q}{n,n+p+q} ^2 I_{pqs}(n-k,0) . \label{mC}
\end{align}
For $k<n$, we~can continue with
\begin{align*}
&= \GG{n+p,n+q}{n,n+p+q} ^2 \Big/ \GG{n-k+p,n-k+q}{n-k,n-k+p+q} ^2 c_{pq,0}(s) \quad\text{by \eqref{mA}} \\
&= \GG{n+p,n+q,n-k,n-k+p+q}{n,n+p+q,n-k+p,n-k+q} ^2 c_{pq,0}(s) \\
&\asymp \frac{(1+p)^{2k} (1+q)^{2k}} {(1+p+q)^{2k}} c_{pq,0}(s),
\end{align*}
as asserted.

For $k=n$,  we have from Theorem~\ref{PE}
\begin{multline*}
 I_{pqs}(n,0) = \GG{n+p+q}{n+p,n+q}^2 \GG{s+1,s+1+2n}{2s+2n+2} (p)_n (q)_n \\
 \int_0^1 \int_0^1 x^{p-1} y^{q-1} (1-x)^{n+s} (1-y)^{n+s}
  \FF21{s+1,n+s+1}{2n+2s+2}{1-xy} \, dx \,dy .
\end{multline*}
For fixed $p,q\ge1$ and $s>-1$, both sides are holomorphic in $n$ in the half-plane $\Re s>\max(-\frac{s+1}2,-p-q)$,
hence the formula remains in force also for $n=0$:
$$ I_{pqs}(0,0) = \GG{p+q}{p,q}^2 \GG{s+1,s+1}{2s+2}
 \int_0^1 \int_0^1 x^{p-1} y^{q-1} (1-x)^s (1-y)^s
  \FF21{s+1,s+1}{2s+2}{1-xy} \, dx \,dy . $$
The last ${}_2\!F_1$ is increasing on the interval~$[0,1]$, with value 1 at~$z=0$, and is known to have
logarithmic growth at~1 (see~\cite[formula~(12) in~\S2.10]{BE}), hence it is $\asymp1-\log(1-z)$ for $0<z<1$. This gives
$$ I_{pqs}(0,0) \asymp \GG{p+q}{p,q}^2 \int_0^1 \int_0^1 x^{p-1} y^{q-1} (1-x)^s (1-y)^s (1-\log xy) \, dx \,dy  $$
(uniformly for all $p,q\ge1$, for any fixed $s>-1$).
Performing the integration yields
$$ I_{pqs}(0,0) \asymp \GG{p+q}{p,q}^2 [ f(p)f(q) - f'(p)f(q) - f(p)f'(q) ],  $$
where we have temporarily denoted, for $x\ge1$,
$$ f(x) := \frac{\Gamma(x)}{\Gamma(x+s+1)}. $$
Note that
$$ -\frac{f'(x)}{f(x)} = -\psi(x)+\psi(x+s+1) \sim (s+1)/x \quad\text{as } x\to+\infty  $$
is a positive function on $x\ge1$ which vanishes as $x\to+\infty$, thanks to the properties
of the digamma function~$\psi=\Gamma'/\Gamma$. Consequently,
$$ 1-\frac{f'(p)}{f(p)} - \frac{f'(q)}{f(q)} \asymp 1 \qquad\text{on }p,q\ge1 ,  $$
and
$$ I_{pqs}(0,0) \asymp \GG{p+q}{p,q}^2 f(p) f(q) = \GG{p+q,p+q}{p,q,p+s+1,q+s+1}. $$
Thus we can continue \eqref{mC} with
\begin{align*}
\frac{c_{pq,n}(s)}{(p+q+1)^{2n}} & \asymp \GG{n+p,n+q}{n,n+p+q}^2 \GG{p+q,p+q}{p,q,p+s+1,q+s+1} \\
&\asymp (1+p)^{2n-s-1} (1+q)^{2n-s-1} (1+p+q)^{-2n} \quad\text{by Stirling's formula} \\
&\asymp \frac{(1+p)^{2n}(1+q)^{2n}}{(1+p+q)^{2n}} c_{pq,0}(s),
\end{align*}
proving the claim.
\end{proof}

\begin{proof}[Proof of Proposition~\ref{PG}] The~claim \eqref{WD} is now immediate from \eqref{TG}, \eqref{ccc} and the last corollary.
\end{proof}

\begin{corollary} \label{PH}
For $n\ge2$, $s>-1$ and $m=0,1,2,\dots,n$, the squared norm in $W^m(\Bn,d\mu_s)$
of $f=\sum_{p,q}f_{pq}$, $f_{pq}\in\bhpq$, \Mh on~$\Bn$, is~equivalent~to
$$ \sum_{p,q} [(p+1)(q+1)]^{2m-s-1} \|f_{pq}\|^2_\pBn . $$
\end{corollary}

\begin{proof} By~Proposition~\ref{PG} and part~(b) of Proposition~\ref{PF},
$$ \|f\|^2_{W^m(\Bn,d\mu_s)} \asymp \sum_{p,q} d_{pqm} \|f_{pq}\|^2_\pBn  $$
with
$$ d_{00,m} = c_{00,0}(s) \asymp1  $$
and
\begin{align*}
d_{pqm} &= \sum_{l=0}^m \sum_{k=0}^l (p+q)^{2(l-k)} [(p+1)(q+1)]^{2k-s-1} \\
&\asymp \sum_{l=0}^m \frac{[(p+q)^2+(p+1)^2(q+1)^2]^l}{(p+1)^{s+1}(q+1)^{s+1}}  \\
&\asymp \frac{[(p+q)^2+(p+1)^2(q+1)^2]^{m}}{(p+1)^{s+1}(q+1)^{s+1}}  \\
&\asymp [(p+1)(q+1)]^{2m-s-1},
\end{align*}
as~asserted.
\end{proof}

\begin{remark} The~last proof shows that in parts (a), (b) of Proposition~\ref{PF}
it~is possible to replace $\sum_{l=0}^m$ by $\sum_{l\in\{0,m\}}$.  \qed
\end{remark}

Let $W^t_\mh(\Bn)$ denote the subspace of all \Mh functions in~$W^t(\Bn)$.

\begin{theorem} \label{PB}
For $0\le t\le n$, $W^t_\mh(\Bn)=\wt\MM_{-2t}$, with equivalent norms.
\end{theorem}

\begin{proof} By~\eqref{VD}, \eqref{VE} and the last corollary with $s=0$,
the~assertion holds for $t=0,1,2,\dots,n$. Since the spaces $\MM_{\#,s}$, $s\in\RR$,
form an interpolation scale (being the completions of the domains of powers $A^{-s/2}$
of the operator $A:\sum_{p,q}f_{pq}\mapsto\sum_{p,q}(p+1)(q+1)f_{pq}$ on~$\MM_0$,
see~\cite[Section~2.1]{LM}) and so do the Sobolev spaces~$W^t$ and, hence, also their 
subspaces~$W^t_\mh$ (see~\cite[Theorem~1 in~\S1.17.1]{Tr}), the result follows by interpolation.
\end{proof}

We~remark that the last theorem in general fails for $t>n$: in~fact, $\bhpq\subset W^t(\Bn)$
for $t\ge n+1$ only if $pq=0$, i.e.~for $t\ge n+1$ the only \Mh functions in $W^t(\Bn)$ are
the pluriharmonic ones. (In~view of known facts like Proposition~1.4 in Graham~\cite{Gr},
this comes as no~surprise.) Indeed, by~\eqref{UC}, any nonzero $f\in\bhpq$ is of the form
$f(r\zeta)=r^{p+q}S_{pq}(r^2)f(\zeta)$, hence $\cN^{n+1} f=\cN^{n+1}[r^{p+q}S_{pq}(r^2)]f(\zeta)$;
but
$$ \cN^{n+1}[r^{p+q}\FF21{p,q}{p+q+n}{r^2}] \approx \frac{2^n\Gamma(p+q+n)}{\Gamma(p)\Gamma(q)} \frac1{1-r} \qquad\text{as }r\nearrow1 $$
fails to be square-integrable near $r=1$ for $pq>0$. We~have not tried to see whether
Theorem~\ref{PB} extends also to some $t$ between $n$ and $n+1$, or to~$t<0$.

Likewise, we made no effort to extend Theorem~\ref{PB} to the weighted Sobolev spaces $W^t(\Bn,d\mu_s)$
with arbitrary $s>-1$. The~argument in the last paragraph shows that, quite generally,
for $s>-1$, $k\in\NN$ and $pq>0$,
$$ c_{pq,n+k}(s) < +\infty \iff s-2k>-1, $$
implying that
$$ k\ge\frac{s+1}2 \implies W^{n+k}_\mh(\Bn,d\mu_s) = W^{n+k}_\ph(\Bn,d\mu_s),  $$
the subspace of all pluriharmonic functions in~$W^{n+k}(\Bn,d\mu_s)$. This suggests that
$W^t_\mh(\Bn,d\mu_s)=W^t_\ph(\Bn,d\mu_s)\iff t\ge n+\frac{s+1}2$.
(Note that, quite generally, it~follows from the Schur lemma applied to the $\Un$ action \eqref{UE}
that
$$ W^t_\mh(\Bn,d\mu_s)=\bigoplus\{\bhpq: \;\bhpq\cap W^t(\Bn,d\mu_s)\neq\{0\} \}  $$
--- where the direct sum is understood with respect to the $W^t$ norm ---
for any $t\in\RR$ and $s>-1$; see Proposition~1 in~\cite{EY}.)

We~remark that for the pluriharmonic functions, the analogue of Theorem~\ref{PB} holds without
any restriction on the order~$t$. This follows easily from the well-known fact that
$$ W^t_\hol(\Bn) = \AA_{-2t} \qquad \forall t\in\RR  $$
in the holomorphic case.

\end{document}